\documentclass[10pt]{article}
\usepackage{amssymb}
\usepackage{amsfonts}
\usepackage{amsmath}
\usepackage[dvips]{graphicx}
\usepackage{latexsym,amsmath,amsfonts,amscd}
\usepackage{epsfig}
\usepackage{changebar}
\usepackage{fancybox}
\title
{Gradient Estimates for
 the Perfect and Insulated Conductivity Problems with Multiple Inclusions}
\author{Ellen ShiTing Bao
\footnote{School of Mathematics, University of Minnesota, 206
Church St SE, Minneapolis, MN 55455, email:~shbao@math.umn.edu
   } \hspace{1cm}
YanYan Li \footnote{Department of Mathematics, Rutgers University,
110 Frelinghuysen Rd. Piscataway, NJ 08854,
email:~yyli@math.rutgers.edu} \hspace{1cm}
 Biao Yin
 \footnote{Department of Mathematics, University of Connecticut,
196 Auditorium Rd. Storrs, CT 06269, email:~yin@math.uconn.edu}}
\date{}

\begin{document}
\newtheorem{Def}{Definition}[section]
\newtheorem{thm}{Theorem}[section]
\newtheorem{lem}{Lemma}[section]
\newtheorem{rem}{Remark}[section]
\newtheorem{prop}{Proposition}[section]
\newtheorem{cor}{Corollary}[section]
\def\av{{\int \hspace{-2.25ex}-} }
\maketitle
\numberwithin{equation}{section}
\bigskip

\setcounter{section}{-1}
\begin{abstract}
In this paper, we study the perfect and the insulated conductivity
problems with multiple inclusions imbedded in a bounded domain in
$\mathbb{R}^n, n\ge 2$.   For these two extreme cases of the
conductivity problems, the gradients of their solutions may blow
up as two inclusions approach each other.   We establish the
gradient estimates for the perfect conductivity problems and an
upper bound of the gradients for the insulated conductivity
problems in terms of the distances between any two closely spaced
inclusions.

\end{abstract}
%%%%%%%%%%%%%%%%%%%%%%%%%%%%%%%%%%%%%%%%%%%%%%%%%%%%%%%%%%%%%%%%%%%%%%%%
%
%     1. Introduction and Statements of Results
%
%%%%%%%%%%%%%%%%%%%%%%%%%%%%%%%%%%%%%%%%%%%%%%%%%%%%%%%%%%%%%%%%%%%%%%%%

\section{Introduction}
  In this paper, a continuation of  \cite{BLY},
 we
establish  gradient estimates for
 the perfect  conductivity problems in
the presence of multiple closely spaced inclusions
 in a bounded
domain in $\mathbb{R}^n$ $(n\ge2)$.  We also establish an upper
bound of the gradients for
 the insulated conductivity problems.
For these two extreme cases of the conductivity problems, the
electric field, which is represented by the gradient of the
solutions, may blow up as the inclusions approach to each other,
the blow-up rates of the electric field have been studied in
\cite{AKLLL,AKL,BLY,Yu,Yu2}. In particular, when there are only
two strictly convex inclusions, and let $\varepsilon$ be the
distance between the two inclusions, then for the perfect
conductivity problem,
 the optimal blow-up rates for the gradients,
 as $\varepsilon$ approaches to zero,
  were established to be $\varepsilon^{-1/2}$,
$(\varepsilon|\ln{\varepsilon}|)^{-1}$ and $\varepsilon^{-1}$ for
$n=2,~3$ and $n\geq 4$ respectively.  A  criteria, in terms of a
functional of boundary data, for the situation where blow-up rate
is realized was also given. See e.g. the introductions of
\cite{BLY} and \cite{Yu2}
for a more detailed description of these results. More
recently, Lim and Yun in \cite{LY} have obtained further estimates with
explicit dependence of the blow-up rates on the size of the
inclusions for the perfect conductivity problem (see also
\cite{AKLLL} for results of this type),  and
 H. Ammari, H. Kang,  H. Lee,
 M. Lim and H. Zribi in \cite{AKLLZ} have given
more
 refined estimates of the gradient of solutions.

 The partial differential equations for the conductivity
problems arise also in the study of composite materials. In
$\mathbb{R}^2$, as explained in \cite{LV}, if we use the bounded
domain to represent the cross-section of a fiber-reinforced
composite and use the inclusions to represent the cross-sections
of the embedded fibers, then by a standard anti-plane shear model,
the conductivity equations can be derived, in which the electric
potential corresponds to the out-of-plane elastic displacement and
the electric field corresponds to the stress tensor.  Therefore,
the gradient estimates for the conductivity problems provide
valuable information about the stress intensity inside the
composite materials.

When conductivities of the inclusions are away from zero and
infinity, the boundedness of the gradients were observed
numerically
 by Babuska, Anderson,
Smith and Levin \cite{BASL}.   Bonnetier and Vogelius \cite{BV}
proved it when the inclusions are two touching balls in
$\mathbb{R}^2$. General results were established
 by Li and Vogelius \cite{LV} for second
order divergence form elliptic equations with piecewise smooth
coefficients, and then by  Li and Nirenberg \cite{LN} for second
order
 divergence form elliptic systems,
including linear system of elasticity,
  with piecewise smooth coefficients.
See also \cite{KSW} and \cite{MOV} for related studies.

\textbf{Acknowledgment:}\hspace{.1cm}We would like to thank Haim
Brezis, Luis Caffarelli, Hyeonbae Kang and Micheal Vogelius for
their suggestions, comments and encouragements to our work.
The work of Y.Y. Li is partially supported by
NSF grant DMS-0701545.

\section{Mathematical set-up and the main results}
Let $\Omega$ be a domain in $\mathbb{R}^n$ with $C^{2,\alpha}$
boundary, $n\geq 2$, $0<\alpha <1$. Let $\{D_{i}\} ~(1\le i\le m)$
be $m$ strictly convex open subsets in $\Omega$ with
$C^{2,\alpha}$ boundaries, $m\ge 2$, satisfying
\begin{equation} \label{domain conditions m inclusions}
\begin{split}
&\text{the principal
 curvature of} ~\partial D_i\geq \kappa_0,\\
&\varepsilon_{ij}:=\text{dist}(D_i, D_j)>0, ~~(i\neq
j)\\
&\text{dist} (D_i,
\partial\Omega)>r_0,~~ \text{diam}(\Omega)<\frac 1{r_0},
 \end{split}
\end{equation}
where $\kappa_0, r_0>0$ are universal constants independent of
$\{\varepsilon_{ij}\}$.  We also assume that the $C^{2,\alpha}$
norms of $\partial D_i$ are bounded by some universal constant
independent of $\{\varepsilon_{ij}\}$. This implies that each
$D_i$ contains a ball of radius  $r^*_0$ for some universal
constant $r^*_0>0$ independent of $\{\varepsilon_{ij}\}$.

We state more precisely what it means by saying that the boundary
of a domain, say $\Omega$,  is $C^{2, \alpha}$ for $0<\alpha<1$:
In a neighborhood of every point of $\partial \Omega$, $\partial
\Omega$ is the graph of some  $C^{2, \alpha}$ function of $n-1$
variables. We define the $C^{2,\alpha}$ norm of $\partial \Omega$,
denoted by $\|\partial \Omega\|_{C^{2,\alpha}}$,
 as the smallest positive number $\frac 1a$ such that in the $2a-$neighborhood
of every point of $\partial \Omega$, identified as $0$ after a
possible translation and rotation of the coordinates so that
$x_n=0$ is the tangent  to $\partial \Omega$ at $0$, $\partial
\Omega$ is given by the graph of a $C^{2,\alpha}$ function,
denoted as $f$, which is defined as $|x'|<a$, the $a-$neighborhood
of $0$ in the tangent plane. Moreover,
$\|f\|_{C^{2,\alpha}(|x'|<a)} \le \frac 1a$.

Denote
$$
\widetilde\Omega :=\Omega\backslash\overline{\cup_{i=1}^m D_{i}}.
$$
Given $\varphi\in C^{1,\alpha}(\partial\Omega)$, the conductivity
problem can be modeled by the following equation:
\begin{equation} \label{eq:k finite}
\left\{ \begin{aligned}
 &div(a_k(x)\nabla u_k)=0 \hspace{0.8cm}\text{in}~\Omega, \\
           &u_k=\varphi
           \hspace{2.6cm}\text{on}~\partial\Omega,
                 \end{aligned}
\right.
\end{equation}
 where $k=(k_1,\ldots,k_m)$ and
\begin{equation}
            a_k(x)=\left\{ \begin{aligned}
           &k_i\in (0, \infty)\hspace{0.6cm}\text{in}~D_{i}, \\
           &1\hspace{2cm}\text{in}~\widetilde\Omega.
          \end{aligned}
\right. \label{0.3}
\end{equation}
The existence and uniqueness of solutions to the above equation is
well known.  Moreover, we have $\|u_k\|_{H^1(\Omega)}\leq
C\|\varphi\|_{C^{1,\alpha}(\partial \Omega)}$ for some constant
$C$ independent of $k$.  Therefore, by passing to a subsequence,
we have $u_k\rightharpoonup u_{\infty}$ in $H^1(\Omega)$ as
$k_i\rightarrow \infty$ for all $1\le i\le m$, where $u_\infty\in
H^1(\Omega)$ is the solution to the following perfect conductivity
problem,
\begin{equation} \label{eq:k +infty m inclusions}
\left\{ \begin{aligned}
           &\Delta u=0 \hspace{2.28cm}in\hspace{0.3cm}\widetilde\Omega, \\
           &u|_{+}=u|_{-}\hspace{1.85cm}on\hspace{0.2cm}\partial D_i,  ~(i=1,2,\ldots,m), \\
           &\nabla u\equiv 0\hspace{2.22cm}in\hspace{0.2cm} D_i ~(i=1,2,\ldots,m), \\
           &\int_{\partial D_i}\frac{\partial u}{\partial\nu}\Big|_{+}=0
           \hspace{0.94cm}~(i=1,2,\ldots,m), \\
           &u=\varphi\hspace{2.49cm}on\hspace{0.2cm}\partial\Omega,
          \end{aligned}
\right.
%\nonumber
\end{equation}
where
$$
\frac{\partial u}{\partial\nu}\Big|_{+}:=\lim_{t\rightarrow
0^{+}}\frac{u(x+t\nu)-u(x)}{t}.
$$
Here and throughout this paper $\nu$ is the outward unit normal to
the domain and the subscript $\pm$ indicates the limit from
outside and inside the domain, respectively.  For the derivation
of the above equation, readers can refer to the Appendix of
$\cite{BLY}$.
Note that the proof there is for $k_1=k_2=\cdots=k_m$, but it works also for the general case with modification.\\
 Since the high stress concentration only occurs in the narrow regions between the
fibers, we only need to focus on those narrow regions.

For $i\neq j$, denote
$$
\text{dist}(x^i_{ij}, x^j_{ij})=\text{dist}(D_i,
D_j)=\varepsilon_{ij}>0, ~x^i_{ij}\in\partial D_i,
~x^j_{ij}\in\partial D_j,$$ and
$$x^0_{ij}:=\frac12(x^i_{ij}+x^j_{ij}).
$$
It is easy to see that there exists some positive constant
$\delta<\frac{1}{4}$ which depends only
 on $\kappa_0,~r_0$ and $\{\|\partial
D_i\|_{C^{2,\alpha}}\}$, but is independent of
$\{\varepsilon_{ij}\}$ such that
\begin{equation}\label{condition: Bij}
\mbox{if}\ \varepsilon_{ij}<2\delta,~ B(x^0_{ij}, 2\delta) \
\mbox{only intersects with}\ D_i\ \mbox{and}\ D_j.
\end{equation}
Denote
\begin{equation}
\rho_n(\varepsilon)=\left\{
\begin{aligned}
           &\frac{1}{\sqrt{\varepsilon}}~~~~~~~~~~~for \hspace{0.2cm} n=2,\\
           &\frac{1}{\varepsilon|\ln{\varepsilon}|}~~~~~~~for \hspace{0.2cm} n=3,\\
           &\frac{1}{\varepsilon}~~~~~~~~~~~~~for \hspace{0.2cm} n \geq4.
\end{aligned}
\right.
\end{equation}

Then we have the following gradient estimates for the perfect conductivity problem\\
\begin{thm} \label{thm:upbdd m inclusions}
Let $\Omega, \{D_i\}\subset \mathbb{R}^n$, $n\ge 2$,
$\{\varepsilon_{ij}\}$ be defined as in (\ref{domain conditions m
inclusions}), $\varphi\in L^{\infty}(\partial\Omega)$, $\delta$ be
the universal constant satisfying $(\ref{condition: Bij})$.
Suppose $u_{\infty}\in H^1(\Omega)$ is the solution to equation
(\ref{eq:k +infty m inclusions}), then for any
$\varepsilon_{ij}<\delta$, we have
 $$\|\nabla
u_\infty\|_{L^{\infty}(\widetilde{\Omega}\cap
B(x^0_{ij},\delta))}\leq
C\rho_n(\varepsilon_{ij})\|\varphi\|_{L^{\infty}(\partial
\Omega)}$$ where $C$ is a constant depending only on $n$,
$\kappa_0$, $r_0$, $\{\|\partial D_i\|_{C^{2,\alpha}}\}$, but
independent of $\varepsilon_{ij}$.
\end{thm}

Note that if $\varepsilon_{ij}\ge\delta$, by the maximum principle
and the boundary estimates of harmonic functions, we immediately
get $\|\nabla u_\infty\|_{L^{\infty}(\widetilde{\Omega}\cap
B(x^0_{ij},\delta) )}\leq
C\|u\|_{L^{\infty}(\widetilde\Omega)}\leq
C\|\varphi\|_{L^{\infty}(\partial \Omega)}$.  Here we have used
the fact that $u_\infty$ is constant on each $\partial D_i$. Then
by Theorem \ref{thm:upbdd m inclusions} and standard boundary
Schauder estimates, see e.g. Theorem 8.33 in \cite{GT}, we have
the global gradient estimates of $u_{\infty}$ in
$\widetilde{\Omega}$.

\begin{cor} \label{cor:upbdd m inclusions}
Let $\Omega, \{D_i\} \subset \mathbb{R}^n$, $n\ge 2$,
$\{\varepsilon_{ij}\}$ be defined as in (\ref{domain conditions m
inclusions}), $\displaystyle\varepsilon:=\min_{i\neq
j}\varepsilon_{ij}>0$, and $\varphi\in
C^{1,\alpha}(\partial\Omega)$, $0< \alpha<1$, and let
$u_{\infty}\in H^1(\Omega)$ be the solution to equation (\ref{eq:k
+infty m inclusions}). Then
$$\|\nabla
u_\infty\|_{L^{\infty}(\widetilde{\Omega})}\leq
C\rho_n(\varepsilon)\|\varphi\|_{C^{1,\alpha}(\partial \Omega)}.$$
where
 $C$ is a constant depending only
 on $n$, $m$, $\kappa_0$,
$r_0$,  $\|\partial \Omega\|_{C^{2,\alpha}}$, $\{\|\partial
D_i\|_{C^{2,\alpha}}\}$, but independent of $\varepsilon$.
\end{cor}
\begin{rem}\label{rem1}
The proof of Theorem \ref{cor:upbdd m inclusions} does not need
$D_i$ and $D_j$ to be strictly convex, the strict convexity is
only used in a fixed neighborhood of $x_{ij}^0$ (The size of the
neighborhood is independent of $\{\varepsilon_{ij}\}$).  In fact,
our proofs of Theorem \ref{thm:upbdd m inclusions} also apply,
with minor modification, to more general situations where
 two closely spaced inclusions, $D_i$ and $D_j$,
 are not necessarily convex near points
on the boundaries where minimal distance $\varepsilon$ is
realized; see discussions after the proof of Theorem
\ref{thm:upbdd m inclusions} in Section 2.
\end{rem}

Next, we study the insulated conductivity problem. Similar to the
perfect conductivity problem, the solution to the insulated
conductivity problem is also the weak limit of $u_k$ in
$H^1(\widetilde \Omega)$ as $k$ approaches to 0. Here we
consider the insulated conductivity problem with anisotropic conductivity.\\
 Let $\Omega,D_i\subset \mathbb{R}^n$, $\varepsilon_{ij}$ be defined as in
(\ref{domain conditions m inclusions}), $\varphi\in
C^{1,\alpha}(\partial\Omega)$, suppose $A(x):=\big(a^{ij}(x)\big)$
is a symmetric matrix function in $\widetilde{\Omega}$, where
$a^{ij}(x)\in C^\alpha(\overline{\widetilde{\Omega}})$ and for
some constants $\Lambda\geq\lambda>0$,
$$\|a^{ij}\|_{C^{\alpha}(\overline{\widetilde{\Omega}})}\leq \Lambda,
\hspace{1cm} a^{ij}(x)\xi_i\xi_j\geq \lambda|\xi|^2, ~~\forall
\xi\in\mathbb{R}^n, x\in\widetilde\Omega.$$ Then the anisotropic
insulated conductivity problem can be described by the following
equation,
\begin{equation} \label{eq:k 0}
\left\{ \begin{aligned}
           &\partial_i(a^{ij}\partial_j u)=0
           \hspace{1.8cm}\text{in}~\widetilde\Omega, \\
           &a^{ij}\partial_j u\nu_i=0
           \hspace{2.15cm}\text{on} ~\partial D_i(i=1,2,\ldots,m), \\
           &u=\varphi\hspace{3.25cm}\text{on}~\partial\Omega.
          \end{aligned}
\right.
\end{equation}
The existence and uniqueness of solutions to equation (\ref{eq:k
0}) are elementary, see the Appendix.

As mentioned before, the blow-up can only occur in the narrow
regions between two closely spaced inclusions. Therefore, we only
derive gradient estimates for the solution to (\ref{eq:k 0}) in
those regions.  Without loss of generality, we consider the
insulated conductivity problem in the narrow region between $D_1$
and $D_2$. Assume
  $$\varepsilon=\text{dist}(D_1, D_2)$$
  After a possible translation and rotation, we may assume
\begin{equation}\label{D1, D2 3}
 ~(\varepsilon/2, 0')\in\partial D_1, ~(-\varepsilon/2, 0')\in \partial D_2.
\nonumber
\end{equation}
Here and throughout this paper by writing $x=(x_1,x')$, we mean
$x'$ is the last $n-1$ coordinates of $x$.

 We denote the narrow region between
 $D_1$ and $D_2$ and its boundary on $\partial D_1$ and $\partial D_2$ as follows
\begin{equation}\label{eq:D1D2}
\begin{aligned}
&\mathcal{O}(r):=\widetilde\Omega\cap \{x\in \mathbb{R}^n\big||x'|<r\}\\
&\Gamma_+:=\partial D_1\cap \{x\in \mathbb{R}^n\big| |x'|<r\}\\
&\Gamma_-:=\partial D_2\cap \{x\in \mathbb{R}^n\big| |x'|<r\}\\
\end{aligned}
\end{equation}
where $r$ is some universal constant depending only on
$\{\|\partial D_i\|_{C^{2,\alpha}}\}$.

With the above notations, we consider the following problem,
\begin{equation} \label{eq:k 0 local harmonic}
\left\{ \begin{aligned}
           &\partial_i(a^{ij}\partial_j u)=0
           ~~~~~~~~~~\text{in}~\mathcal{O}(r), \\
           &a^{ij}\partial_j u\nu_i=0
            ~~~~~~~\text{on} ~ \Gamma_+\cup \Gamma_-.
          \end{aligned}
\right.
\end{equation}

Then we have:
\begin{thm} \label{thm:upbdd k 0}
 If $u_0\in H^1(\mathcal{O}(r))$
is a weak solution of (\ref{eq:k 0 local harmonic}), then
\begin{equation} \label{eq:bd:local}
|\nabla u_0(x)|\leq
\frac{C\|u_0\|_{L^{\infty}(\mathcal{O}(r))}}{\sqrt{\varepsilon+|x'|^2}},~~
\textrm{for all}~~ x\in\mathcal{O}(\frac{r}{2}).
\end{equation}
where $C$ is a constant depending only on $n$, $\kappa_0$, $r_0$,
$\Lambda$, $\lambda$, $r$ and $\|\partial
D_i\|_{C^{2,\alpha}}(i=1,2)$, but independent of $\varepsilon$.
\end{thm}
\begin{rem}
 Theorem \ref{thm:upbdd k 0} also remains true for general second
order elliptic systems, its proof is essentially the same as for
the equations.
\end{rem}
 A consequence of Theorem \ref{thm:upbdd k 0} is the following
global gradient estimates for the insulated conductivity problem.
\begin{cor}\label{cor:upbdd k 0}
Let $\Omega, \{D_i\} \subset \mathbb{R}^n$, $\{\varepsilon_{ij}\}$
be defined as in (\ref{domain conditions m inclusions}),
$\displaystyle\varepsilon:=\min_{i\neq j}\varepsilon_{ij}>0$, and
$\varphi\in C^{1,\alpha}(\partial\Omega)$, let $u_0\in
H^1(\widetilde{\Omega})$ be the weak solution to equation
(\ref{eq:k 0}), then
\begin{equation}
\|\nabla u_0\|_{L^{\infty}(\widetilde\Omega)}\leq
\frac{C}{\sqrt\varepsilon}\|\varphi\|_{C^{1,\alpha}(\partial
\Omega)}.
\end{equation}
where $C$ is a constant depending only on $n$, $\kappa_0$, $r_0$,
$\|\partial \Omega\|_{C^{2,\alpha}}$, $\{\|\partial
D_i\|_{C^{2,\alpha}}\}$, but independent of $\varepsilon$.
\end{cor}
Note that throughout this paper we often use $C$ to denote
different constants, but all these constants are independent of
$\varepsilon$.

 The paper is organized as follows. In Section 2 we consider the
perfect conductivity problem and prove Theorem \ref{thm:upbdd m
inclusions}. In Section 3 we show Theorem \ref{thm:upbdd k 0} for
the insulated case. Finally in the Appendix we present some
elementary results for the insulated conductivity
problem.\\

\section{The perfect conductivity problem with multiple inclusions}

In this section, we consider the perfect conductivity problem
(\ref{eq:k +infty m inclusions}). Note that from equation (\ref{eq:k
+infty m inclusions}), we know that
 $u\equiv C_i~\text{on}~\overline {D}_i, ~1\leq i\leq m,$
 where $\{C_i\}$ are some unknown constants. In order to prove Theorem \ref{thm:upbdd m inclusions}, we first
estimate $|C_i-C_j|$ for $1\leq i\neq j\leq m$, which later will
allow us to control the gradient of $u$ in the narrow region
between $D_i$ and $D_j$.
%%%%%%%%%%%%%%%%%%%%%%%%%%%%%%%%%%%%%%%%%%%%%%%%%%%%%%%%%%%%%%%%%%%%%%%%
%
%     2. Estimate of C_i-C_j
%
%%%%%%%%%%%%%%%%%%%%%%%%%%%%%%%%%%%%%%%%%%%%%%%%%%%%%%%%%%%%%%%%%%%%%%%%
\subsection{A Matrix Result}
 To estimate
$|C_i-C_j|$, the following proposition plays a crucial role.

Let $m$ be a positive integer, $P=(p_{ij})$ an $m \times m$ real
symmetric matrix satisfying,
\begin{equation}
\begin{aligned}
   &(A1)~ p_{ij}=p_{ji}\leq 0 ~(i \neq j);\\
   &(A2)~0<r_1\leq\bar p_i:=\sum_{j=1}^{m}p_{ij}\leq r_2,~~~~~~\\
\end{aligned}
\nonumber
\end{equation}
where $r_1$ and $r_2$ are some positive constants.
\begin{rem}
 An $m\times m$ matrix $P$ satisfying $|p_{ii}|>\sum\limits_{j\neq
i} |p_{ij}|$ is called a diagonally dominant matrix. Such a matrix
is nonsingular, see \cite{GV}. $(A1)$ and $(A2)$ imply that the
matrix $P$ is diagonally dominant.
\end{rem}

\begin{prop} \label{prop:system P}
Let $P=(p_{ij})$ be an $m \times m$ real symmetric matrix
satisfying $(A1)$ and $(A2)$, $m\ge 1$. For $\beta\in
\mathbb{R}^m$, let $\alpha$ be the solution of
\begin{equation} \label{system:P}
P\alpha=\beta,
\end{equation}
then
\begin{equation}
|\alpha_i-\alpha_j|\leq
m(m-1)\frac{r_2}{r_1}\frac{|\beta|}{|p_{ij}|+r_1},
\end{equation}
\end{prop}
where $|\beta|=\max\limits_i |\beta_i|$.

Before proving the proposition, we introduce the following lemmas. \\
Denote
\begin{equation}\label{matrix eo set}
\begin{aligned}
&\mathcal{I}(l)=\{\textrm{all}~ l\times l ~\textrm{diagonal
matrices whose diagonal entries are $1$ or $-1$}\},\\
 &\mathcal{I}_e(l)=\{\overline{I}\in \mathcal{I}(l)\big|\overline{I}~ \textrm{has even numbers of $-1$ in its diagonal}\},\\
&\mathcal{I}_o(l)=\{\overline{I}\in \mathcal{I}(l)\big|\overline{I} ~\textrm{has odd numbers of $-1$ in its diagonal}\}.\\
\end{aligned}
\nonumber
\end{equation}
\begin{lem} \label{lem:matrix identity}
For any $x\in\mathbb{R}$ and any $l\times l$ matrix $A$, $l\ge 1$,
$$\sum_{\overline{I}\in\mathcal{I}_e(l)}\det{(xI+\overline{I}A)}\equiv
2^{l-1}(x^l+\det{A});$$
$$\sum_{\overline{I}\in\mathcal{I}_o(l)}\det{(xI+\overline{I}A)}\equiv
2^{l-1}(x^l-\det{A}).$$
\end{lem}
\emph{Proof:}\hspace{.2cm} We prove it by induction. The above
identities can be easily checked for $l=1$. Suppose that the above
identities stand for $l=k-1\ge 1$, we will prove them for $l=k$.
Observe that the above identities hold when $x=0$. To prove them
for all $x$, it suffices to show that the derivatives with respect
to $x$ in both sides of the identities coincide. Since for any
$\overline{I}\in \mathcal{I}(k)$,
$$(\det{(xI+\overline{I}A)})'=\sum\limits_{i=1}^{k}\det{(x I+\overline{I}_i A_{i})}$$
where $A_{i}$ and $\overline{I}_i$ are the submatrices obtained by
eliminating the ith row and the ith column of $A$ and
$\overline{I}$ respectively. \\
Notice that if $\overline{I}$ runs through all the elements of
$\mathcal{I}_e(k)$, $\overline{I}_i$ will run through all the
elements of $\mathcal{I}(k-1)$ for every fixed
$i\in\{1,2,\ldots,k\}$, so we have
\begin{equation}
\begin{aligned}
&~~\sum_{\overline{I}\in\mathcal{I}_e(k)}(\det{(xI+\overline{I}A)})'\\
&=\sum\limits_{i=1}^{k}\big(\sum_{\overline{I}\in\mathcal{I}_e(k-1)}
\det{(x
I+\overline{I}A_{i})}+\sum_{\overline{I}\in\mathcal{I}_o(k-1)}
\det{(xI+\overline{I} A_{i})}\big) \\
&=\sum\limits_{i=1}^{k}\big(2^{k-2}(x^{k-1}+\det{A_{i}})+2^{k-2}(x^{k-1}-\det{A_{i}})\big)~~~~~~(\textrm{By induction})\\
&=k2^{k-1}x^{k-1}=2^{k-1}(x^k+\det{A})'.
\end{aligned}
\nonumber
\end{equation}
 Therefore, we have proved the first identity.  The second one follows from the first one by changing the sign of one row of $A$.\\
As a consequence of Lemma \ref{lem:matrix identity}, we have
\begin{cor}\label{cor:A l*l}
Let $A$ be an $l\times l$ matrix, if $\det{(I+\overline{I}A)}\ge
0$ for any $\overline{I}\in\mathcal{I}(l)$, then $|\det A|\le 1$.
\end{cor}

\begin{lem}\label{lem:Q m*l}
Given integers $m>l\geq 1$, let $Q=(q_{ij})$ be an $m\times l$
real matrix which satisfies, for $j=1,2,\ldots,l$,
\begin{equation} \label{matrix condition}
q_{jj}>\sum_{i\neq j}|q_{ij}|.
\end{equation}
Let $\mathcal{A}$ be the set of all $l\times l$ submatrices of the
above matrix $Q$ and $S_1\in \mathcal{A}$ the matrix obtained from
the first $l$ rows of $Q$, then we have
$$\det S_1=\max_{S\in \mathcal{A}}|\det S|.$$
\end{lem}
\emph{Proof}: For any $S \in \mathcal{A}$, by rearranging the
order of its rows we do not change $|\det S|$. Thus we can treat
$S$ as a matrix obtained by replacing some rows of $S_1$ by some
other rows of $Q$. Note that $S$ and $S_1$ could have no rows in
common, which means $S$ is obtained by replacing all the rows
of $S_1$ by some other rows of $Q$. \\
Given any $\overline{I}\in\mathcal{I}(l)$, we claim:
$$ \det{(S_1+\overline{I}S)}\ge 0$$
\emph{Proof of the claim:}~There are two cases between $S_1$ and $S$:\\
Case 1.~$S_1$ and $S$ have no rows in common. Then by (\ref{matrix
condition}), we know that $S_1+\overline{I}S$ is diagonally dominant,
therefore $\det{(S_1+\overline{I}S)}>0$.\\
Case 2.~$S_1$ and $S$ have some common rows, denote the order of
these rows by $1\le i_1<\cdots< i_s\le l, 1\le s\le l$. If row
$i_{s_0}$ of $\overline{I}S$ is opposite to row $i_{s_0}$ of $S$
for some $1\le s_0\le s$, then row $i_{s_0}$ of
$S_1+\overline{I}S$ is 0, therefore $\det{(S_1+\overline{I}S)}=0$.
Otherwise row $i_t$ of $\overline{I}S$ is the same as that of $S$
and $S_1$ for any $1\le t\le s$, then we take out the common
factors $2$ in these rows when we compute
$\det{(S_1+\overline{I}S)}$, thus we have
$$\det{(S_1+\overline{I}S)}=2^s\det{(S_1+\overline{I}\hat{S})},$$
where $\hat{S}$ is the matrix obtained by replacing row $i_t$ of
$S$ by 0 for any $1\le t\le s$. We know that
$S_1+\overline{I}\hat{S}$ is diagonally dominant according to
(\ref{matrix condition}), then
$\det{(S_1+\overline{I}\hat{S})}>0$, it yields that
$\det{(S_1+\overline{I}S)}>0$. Therefore, the claim is proved. \\
Since $\det{S_1}>0$ and
$$\det{(S_1+\overline{I}S)}=\det{(I+\overline{I}SS^{-1}_1)}\det{S_1}$$
we have, by the claim, that for any $\overline I\in
\mathcal{I}(l)$,
$$\det{(I+\overline{I}SS^{-1}_1)}\ge 0$$
By Corollary \ref{cor:A l*l}, we have
$$|\det{(S S^{-1}_1)}|\leq 1$$
therefore
$$\det{S_1}\ge |\det{S}|.$$

Now we are ready to prove Proposition \ref{prop:system P}.\\

\emph{Proof of Proposition \ref{prop:system P}:}\hspace{.1cm} For
$m=1$ the inequality is automatically true. For $m=2$, we have, by
Cramer's rule,
\begin{equation}
\alpha_1-\alpha_2=\frac{\left| \begin{array}{cc}
\beta_1 & p_{12}  \\
\beta_2 & p_{22}
\end{array}\right|}{\left|
\begin{array}{cc}
                p_{11} & p_{12}  \\
                p_{21} & p_{22}
\end{array}\right|}- \frac{\left| \begin{array}{cc}
p_{11} & \beta_1   \\
p_{21} & \beta_2
\end{array}\right|}{\left|
\begin{array}{cc}
                p_{11} & p_{12}  \\
                p_{21} & p_{22}
\end{array}\right|}=\frac{\left| \begin{array}{cc}
\beta_1 & \bar p_1  \\
\beta_2 & \bar p_2
\end{array}\right|}{\left|
\begin{array}{cc}
                p_{11} & p_{12}  \\
                p_{21} & p_{22}
\end{array}\right|}\\
\nonumber
\end{equation}
Since $r_1\leq\bar p_i\leq r_2$ by Condition (A2),
$$\left| \begin{array}{cc}
\beta_1 & \bar p_1  \\
\beta_2 & \bar p_2
\end{array}\right|=\beta_1\bar p_2-\beta_2\bar p_1\leq 2r_2|\beta|$$
On the other hand, by Condition (A1) and (A2)
$$ \left|
\begin{array}{cc}
                p_{11} & p_{12}  \\
                p_{21} & p_{22}
\end{array}\right|=\left|
\begin{array}{cc}
                \bar p_1 & p_{12}  \\
                \bar p_2 & p_{22}
\end{array}\right|=\bar p_1 p_{22}
-\bar p_2 p_{12}\ge \bar p_1 p_{22}\ge r_1(r_1+|p_{12}|). $$
Therefore, Proposition \ref{prop:system P} for $m=2$ follows from
the above.

For $m\ge 3$, we only estimate $|\alpha_1-\alpha_2|$ since the
other estimates can be obtained by switching columns of $P$.

Since $\alpha$ satisfies (\ref{system:P}), by Cramer's rule, we
have:
\begin{equation}
\begin{aligned}
\alpha_1-\alpha_2&=\frac{\left| \begin{array}{ccccc}
\beta_1 & p_{12}  & \cdots & p_{1m} \\
\beta_2 & p_{22}  & \cdots & p_{2m} \\
\vdots & \vdots  & \ddots & \vdots \\
\beta_m & p_{m2} & \cdots & p_{mm}
\end{array}\right|}{\left|
\begin{array}{cccc}
                p_{11} & p_{12} & \cdots & p_{1m} \\
                p_{21} & p_{22} & \cdots & p_{2m} \\
                \vdots & \vdots & \ddots & \vdots \\
                p_{m1} & p_{m2} & \cdots & p_{mm}\end{array}
                \right|}-
\frac{\left| \begin{array}{ccccc}
p_{11} & \beta_1  & \cdots & p_{1m} \\
p_{21} & \beta_2  & \cdots & p_{2m} \\
\vdots & \vdots  & \ddots & \vdots \\
p_{m1} & \beta_m & \cdots & p_{mm}
\end{array}\right|}{\left|
\begin{array}{cccc}
                p_{11} & p_{12} & \cdots & p_{1m} \\
                p_{21} & p_{22} & \cdots & p_{2m} \\
                \vdots & \vdots & \ddots & \vdots \\
                p_{m1} & p_{m2} & \cdots & p_{mm}\end{array}
                \right|}\\
&=\frac{\left| \begin{array}{ccccc}
\beta_1 & p_{11}+p_{12} & p_{13} & \cdots & p_{1m} \\
\beta_2 & p_{21}+p_{22} & p_{23} & \cdots & p_{2m} \\
\beta_3 & p_{31}+p_{32} & p_{33} & \cdots & p_{3m} \\
\vdots & \vdots & \vdots & \ddots & \vdots \\
\beta_m & p_{m1}+p_{m2} & p_{m3} & \cdots & p_{mm}
\end{array}\right|}{\left|
\begin{array}{cccc}
                p_{11} & p_{12} & \cdots & p_{1m} \\
                p_{21} & p_{22} & \cdots & p_{2m} \\
                \vdots & \vdots & \ddots & \vdots \\
                p_{m1} & p_{m2} & \cdots & p_{mm}\end{array} \right|}
\end{aligned}
\nonumber
\end{equation}
By adding the last $(m-2)$ columns of the matrix in the numerator
to its second column, we have
\begin{equation}
\alpha_1-\alpha_2=\frac{\left| \begin{array}{ccccc}
\beta_1 & \bar p_{1} & p_{13} & \cdots & p_{1s} \\
\beta_2 & \bar p_{2} & p_{23} & \cdots & p_{2s} \\
\beta_3 & \bar p_{3} & p_{33} & \cdots & p_{3s} \\
\vdots & \vdots & \vdots & \ddots & \vdots \\
\beta_m & \bar p_{m} & p_{m3} & \cdots & p_{mm}
\end{array}\right|}{
\left|
\begin{array}{cccc}
                p_{11} & p_{12} & \cdots & p_{1m} \\
                p_{21} & p_{22} & \cdots & p_{2m} \\
                \vdots & \vdots & \ddots & \vdots \\
                p_{m1} & p_{m2} & \cdots & p_{mm}\end{array}
                \right|}:=\frac{\det \widetilde{P}}{\det P}.
\nonumber
\end{equation}

Next we estimate the determinants of the above two matrices
separately.

Expanding $\det P$ with respect to the first column, we have
$$\det P=\sum_{j=1}^m p_{j1}P_{j1}$$
where $P_{ji}$ is the cofactor of $p_{j1}$.\\
 Applying Lemma
\ref{lem:Q m*l} to the $m\times (m-1)$ matrix obtained by
eliminating the first column of $P$, we know that, among the
cofactors $P_{j1}$, $P_{11}>0$ has the largest absolute value.
Since $p_{j1}=p_{1j}\le0~(j\neq 1)$ and $p_{11}>0$ by condition
$(A1)$ and $(A2)$, we have
\begin{equation}
\det P\ge \sum_{j=1}^m p_{j1}P_{11}=\bar p_1 P_{11}.
 \nonumber
\end{equation}
For the same reason, we have
\begin{equation}
P_{11}=\left| \begin{array}{ccc}
                p_{22} & \cdots & p_{2m} \\
                \vdots & \ddots & \vdots \\
                p_{m2} & \cdots & p_{mm} \\
                \end{array} \right| \geq\big(\sum_{j=2}^m p_{2j}\big)\left| \begin{array}{ccc}
                p_{33} & \cdots & p_{3m} \\
                \vdots & \ddots & \vdots \\
                p_{m3} & \cdots & p_{mm} \\
                \end{array} \right|.
                \nonumber
\end{equation}
Combining the above two inequalities and using condition $(A1)$
and $(A2)$, we have

\begin{equation}
\label{eqn:P}
\begin{aligned}
\det P&\geq\bar p_1\sum_{j=2}^m p_{2j}\left| \begin{array}{ccc}
                p_{33} & \cdots & p_{3m} \\
                \vdots & \ddots & \vdots \\
                p_{m3} & \cdots & p_{mm} \\
                \end{array} \right|
                =\bar p_1(\bar p_2-p_{21})\left| \begin{array}{ccc}
                p_{33} & \cdots & p_{3m} \\
                \vdots & \ddots & \vdots \\
                p_{m3} & \cdots & p_{mm} \\
                \end{array} \right|\\
           &\geq r_1(|p_{12}|+r_1)\left| \begin{array}{ccc}
                p_{33} & \cdots & p_{3m} \\
                \vdots & \ddots & \vdots \\
                p_{m3} & \cdots & p_{mm} \\
                \end{array} \right|.
\end{aligned}
\end{equation}

 By Laplace expansion, see e.g. page 130 of \cite{V}, we can
 expand $\det\widetilde{P}$ with respect to the first two columns of
 $P$, namely,
 \begin{equation}\label{eq:widetilde P}
 \det\widetilde{P}=\sum_{i_1,i_2}\left| \begin{array}{cc}
                \beta_{i_1} &\bar p_{i_1} \\
                \beta_{i_2} & \bar p_{i_2} \\
                \end{array} \right|\widetilde P_{i_1 i_2 12},
                \end{equation}
 where $1\le i_1<i_2\le m$ and $\widetilde P_{i_1 i_2 1 2}$ is the cofactor
 of the 2nd-order minor in row $i_1, i_2$ and column $1,2$ of $\widetilde
 P$.\\
  Applying Lemma \ref{lem:Q m*l} to the $m\times (m-2)$ matrix obtained by eliminating the
first $2$ columns of $\widetilde{P}$, we know that, among all
those cofactors,
$$\left| \begin{array}{ccc}
                p_{33} & \cdots & p_{3m} \\
                \vdots & \ddots & \vdots \\
                p_{m3} & \cdots & p_{mm} \\
                \end{array} \right|$$
has the largest absolute value. Since $0<\bar p_i\leq r_2$ by
condition $(A2)$,
$$\left| \begin{array}{cc}
                \beta_{i_1} &\bar p_{i_1} \\
                \beta_{i_2} & \bar p_{i_2} \\
                \end{array} \right|\le 2r_2|\beta|,$$
then by (\ref{eq:widetilde P}), we have
\begin{equation}\label{eqn:P tilde}
\big|\det \widetilde{P}\big|\leq m(m-1)r_2|\beta|\left|
\begin{array}{ccc}
                p_{33} & \cdots & p_{3m} \\
                \vdots & \ddots & \vdots \\
                p_{m3} & \cdots & p_{mm} \\
                \end{array} \right|.
\end{equation}

By (\ref{eqn:P}) and (\ref{eqn:P tilde}), we have
$$
{|\alpha_1-\alpha_2|=\frac{|\det \widetilde{P}|}{|\det P|}\leq
m(m-1)\frac{r_2}{r_1}\frac{|\beta|}{|p_{12}|+r_1}.}$$

\subsection{Proof of Theorem \ref{thm:upbdd m inclusions}}
As in \cite{BLY}, we decompose $u_\infty$ into $m+1$ parts:
\begin{equation} \label{rep:u m inclusions}
u_\infty=v_{0}+\sum_{i=1}^m C_i v_i,
\end{equation}
where $v_i\in H^1(\widetilde{\Omega})~(i=0,1,2,\ldots,m)$ are
determined by the following equations:\\
for $i=0$,
\begin{equation} \label{eq:v_0 m inclusions}
\left\{ \begin{aligned}
          \Delta v_0&=0 \hspace{1.8cm}in\hspace{0.3cm}\widetilde{\Omega}, \\
           v_0&=0\hspace{1.8cm}on\hspace{0.2cm}\partial D_1,~\partial D_2,\ldots~\partial D_m,\\
           v_0&=\varphi\hspace{1.74cm}on\hspace{0.2cm}\partial\Omega.
          \end{aligned}
\right.
%\nonumber
\end{equation}
for $i=1,2,\ldots, m$,
\begin{equation} \label{eq:v_i m inclusions}
\left\{ \begin{aligned}
           \Delta v_i&=0 \hspace{2cm}in\hspace{0.3cm}\widetilde{\Omega}, \\
           v_i&=1\hspace{2cm}on\hspace{0.2cm}\partial D_i,\\
           v_i&=0\hspace{2cm}on\hspace{0.2cm}\partial D_j, ~\textrm{for}~ j\neq i,  \\
           v_i&=0\hspace{2cm}on\hspace{0.2cm}\partial\Omega.
          \end{aligned}
\right.
%\nonumber
\end{equation}
Since $u_\infty$ satisfies the integral conditions in equation
(\ref{eq:k +infty m inclusions}), using the decomposition formula
(\ref{rep:u m inclusions}), we know that the vector $(C_1, C_2,
\ldots, C_m)$ satisfies the following system of linear equations
\begin{equation} \label{system:C_1 C_2 C_m inclusions}
\begin{pmatrix} a_{11} & a_{12} & \cdots & a_{1m} \\
                a_{21} & a_{22} & \cdots & a_{2m} \\
                \vdots & \vdots & \ddots & \vdots \\
                a_{m1} & a_{m2} & \cdots & a_{mm} \\
                \end{pmatrix}
\begin{pmatrix} C_{1} \\
                C_{2} \\
                \vdots \\
                C_{m} \\
                \end{pmatrix}
=
\begin{pmatrix} b_{1} \\
                b_{2} \\
                \vdots \\
                b_{m} \\
                \end{pmatrix}
%\nonumber
\end{equation}
where
\begin{align}
\label{def:a_ij m inclusions} &a_{ij}:=\int_{\partial
D_j}\frac{\partial
v_i}{\partial\nu}, ~~(i,j=1,2,\ldots, m),\\
\label{def:b_i m inclusions} &b_i:=-\int_{\partial
D_i}\frac{\partial v_0}{\partial\nu}, ~~~~(i=1,2,\ldots, m).
\end{align}

Similar to the two inclusions case in \cite{BLY}, we first
investigate the properties of $v_i~(i=0,1,\cdots,m)$, the matrix
$A=(a_{ij})$ and the vector $b$ defined by (\ref{def:a_ij m
inclusions}) and (\ref{def:b_i m inclusions}). Here we state the
following lemma, for its proof, readers may refer to Lemma $2.4$
in \cite{BLY}.

\begin{lem}\label{lm:ab m inclusions}
For $1\leq i,~ j\leq m$, let $a_{ij}$ and $b_{i}$ be defined by
(\ref{def:a_ij m inclusions}) and (\ref{def:b_i m inclusions}),
then they satisfy the following:
\begin{enumerate}
   \item[(1)]   $~a_{ii}<0 , ~~a_{ij}=a_{ji}>0 ~(i\neq j)$,
   \item[(2)]   $\displaystyle~-C\leq \sum_{1\leq j\leq m} a_{ij}\leq -\frac{1}{C}$,
   \item[(3)]   $~|b_i|\leq C\|\varphi\|_{L^{\infty}(\partial \Omega)}$,
\end{enumerate}
where $C>0$ is a universal constant depending only on $n$,
$\kappa_0$, $r_0$, $\|\partial \Omega\|_{C^{2,\alpha}}$, but
independent of $\varepsilon_{ij}$.
\end{lem}

\begin{rem}\label{rem1 3}
From property (1) and (2) in Lemma \ref{lm:ab m inclusions}, we
know that $A$ is diagonally dominant, therefore it is nonsingular.
\end{rem}

\begin{lem}\label{lm grad v_0,i m inclusions}
Let $v_0, v_i (i=1,\ldots,m)$ be the solutions of equations
(\ref{eq:v_0 m inclusions}) and ({\ref{eq:v_i m inclusions}})
respectively, $\delta$ is the constant satisfying (\ref{condition:
Bij}), then there exists a universal constant $C$ depending only
on $n$, $m$, $r_0$, $\kappa_0$, $\|\partial D_i\|_{C^{2,\alpha}}$
and $\|\partial \Omega\|_{C^{2,\alpha}}$, but independent of
$\{\varepsilon_{ij}\}$ such that,
\begin{enumerate}
   \item[(1)]    $\|\nabla v_0\|_{L^{\infty}(\widetilde{\Omega})}\leq
   C$;
   \item[(2)]   $\|\nabla
v_i\|_{L^{\infty}(B(x_{ij}^0,\delta)\cap\widetilde{\Omega})}\leq
\frac{C}{\varepsilon_{ij}}$ ~if $\varepsilon_{ij}<\delta$;
   \item[(3)]   $|\nabla
v_i|\leq C$ ~~on ~~ $\widetilde
\Omega\setminus{\big(\bigcup_{j\neq
i,\varepsilon_{ij}<\delta}B(x_{ij}^0,\delta)\big)}$.
\end{enumerate}
\end{lem}
\emph{Proof}: \hspace{.2cm}The proof of (1) is the same as the
proof of Lemma 2.3 in \cite{BLY}.  Since
$\|v_i\|_{L^\infty(\widetilde{\Omega})}=1$, $\delta$ is the
constant satisfying (\ref{condition: Bij}), then if
$\varepsilon_{ij}<\delta$, then by (\ref{condition: Bij}), we know
that $B(x_{ij}^0,\delta)$ only intersects with $D_i$ and $D_j$,
and $B(x_{ij}^0,\delta)$ is at least $\delta$ away from other
inclusions. Then (2) just follows from the maximum principle and
standard boundary estimates for harmonic functions. For the same
reason, to prove (3), we only need to prove $\|\nabla
v_i\|_{L^\infty(B(x^0_{kl},\delta)\cap\widetilde\Omega)}\leq C$ if
$k,l\neq i$ and $\varepsilon_{kl}<\delta$.  Without loss of
generality, we assume $k=1,l=2, i=3$. Let $\widetilde v_3$ be the
solution of the following equation,
\begin{equation}
\left\{ \begin{aligned}
          \Delta \widetilde v_3&=0 \hspace{1.8cm}in\hspace{0.3cm}\Omega\setminus {\overline {D_1\cup D_3}}, \\
          \widetilde v_3&=0\hspace{1.8cm}on\hspace{0.2cm}\partial D_1,\\
           \widetilde v_3&=1\hspace{1.8cm}on\hspace{0.2cm}\partial D_3,\\
           \widetilde
           v_3&=0\hspace{1.74cm}on\hspace{0.2cm}\partial\Omega.
          \end{aligned}
\right. \nonumber
\end{equation}
Then we have $\widetilde v_3\ge v_3$ on $\partial\widetilde
\Omega$, by the maximum principle, $\widetilde v_3\ge v_3$ in
$\widetilde \Omega$.  Since $\widetilde v_3=v_3=0$ on $\partial
D_1$, we have
$$\frac{\partial \widetilde v_3}{\partial \nu}\ge\frac{\partial
v_3}{\partial\nu}\ge 0.$$
But $|\nabla \widetilde v_3|<C$ on
$\partial D_1\cap B(x_{12}^0,\delta)$ by the boundary estimates of
harmonic functions, then we have
\begin{equation}\label{eq:v3 D1}
\|\nabla  v_3\|_{L^{\infty}(\partial D_1\cap
B(x_{12}^0,\delta))}=\|\frac{\partial v_3}{\partial\nu}
\|_{L^{\infty}(\partial D_1\cap B(x_{12}^0,\delta))}<C.
\end{equation}
Similarly, we have
\begin{equation}\label{eq:v3 D2}
\|\nabla  v_3\|_{L^{\infty}(\partial D_2\cap
B(x_{12}^0,\delta))}=\|\frac{\partial v_3}{\partial\nu}
\|_{L^{\infty}(\partial D_2\cap B(x_{12}^0,\delta))}<C.
\end{equation}

 Furthermore, by gradient estimates and boundary
estimates of harmonic functions, we have
\begin{equation} \label{eq:v3 side}
\|\nabla v_3\|_{L^\infty({\partial B(x_{12}^0,\delta)\cap
\widetilde\Omega})}<C.
\end{equation}
Since $\nabla v_3$ is still harmonic function on $
B(x_{12}^0,\delta)\cap \widetilde\Omega$, by (\ref{eq:v3 D1}),
(\ref{eq:v3 D2}) and (\ref{eq:v3 side}) and the maximum principle,
we have
\begin{equation}
\|\nabla v_3\|_{L^{\infty}(\widetilde{\Omega}\cap
B(x^0_{12},\delta))}<C. \nonumber
\end{equation}

Next, we derive some further estimates of $A=(a_{ij})$.

\begin{lem}\label{order aij}
Let $a_{ij}$ be defined as in (\ref{def:a_ij m inclusions}), then
there exists a universal constant $C>0$, depending only on $n$,
$r_0$, $\kappa_0$, $\|\partial D_i\|_{C^{2,\alpha}}$ and
$\|\partial \Omega\|_{C^{2,\alpha}}$, but independent of
$\{\varepsilon_{ij}\}$, such that for $1\le i\neq j\le m$,
\begin{align*}
-\frac{C}{\sqrt{\displaystyle\min_{k\neq
i}\varepsilon_{ik}}}<a_{ii}<-\frac{1}{C\sqrt
{\displaystyle\min_{k\neq i}\varepsilon_{ik}}},\hspace{1.5cm}
\frac{1}{C\sqrt {\varepsilon_{ij}}}<a_{ij}<\frac{C}{\sqrt
{\varepsilon_{ij}}},~~~~~~&\text{for}~ n=2,\\
-C|\ln(\min_{k\neq
i}\varepsilon_{ik})|<a_{ii}<-\frac1C|\ln(\min_{k\neq
i}\varepsilon_{ik})|,\hspace{0.5cm} \frac1C|\ln
\varepsilon_{ij}|<a_{ij}<C|\ln\varepsilon_{ij}|,~~~~&\text{for}~ n=3,\\
-C<a_{ii}<-\frac1C,
\hspace{3.5cm}\frac1C<a_{ij}<C,~~~~~~~~~~&\text{for}~ n\geq 4.
\end{align*}
\end{lem}
\emph{Proof}:\hspace{.2cm} Without loss of generality, we assume
$i=1, j=2$. The proof of the estimates for $a_{11}$ is the same as
that in Lemma $2.5$, Lemma $2.6$, and Lemma $2.7$ in \cite{BLY}.
Here we prove the estimate for $a_{12}$.  In the following, we use
$C$ to denote some universal constant depending only on $n$,
$r_0$, $\kappa_0$, $\|\partial D_i\|_{C^{2,\alpha}}$ and
$\|\partial \Omega\|_{C^{2,\alpha}}$, but independent of
$\{\varepsilon_{ij}\}$.

 Notice that if
$\varepsilon_{12}$ is larger than some universal constant, then
the proof is trivial. Therefore, we can assume
$\varepsilon_{12}<\delta$, where $\delta<1/4$ is the universal
constant satisfying (\ref{condition: Bij}). By (\ref{condition:
Bij}), we know that $B(x^0_{12},\delta)$ only intersects with
$D_1$ and $D_2$.

Denote \[\Gamma_i:=\partial D_i\cap B(x^0_{12},\delta)~(i=1,2),
~~\Gamma_3:=
\partial
B(x^0_{12},\delta)\setminus{(D_1\cup D_2)}\]

Since $B(x^0_{12},2\delta)$ does not intersect with $D_i(i\ge3)$
or $\partial \Omega$ by (\ref{condition: Bij}),
  then $$dist(\Gamma_3,\cup^m_{i=3}\partial
 D_i)>\delta,~~dist(\Gamma_3,\partial \Omega)>\delta,$$
 by standard gradient estimates
and boundary estimates for harmonic functions,
 we have
\begin{equation}\label{eq:v1 gamma3}
 \|\nabla
v_1\|_{L^{\infty}(\Gamma_3)}<C
\end{equation}
By Lemma \ref{lm grad v_0,i m inclusions}, we have $\|\nabla
 v_1\|_{L^{\infty}(\partial D_2\setminus \Gamma_2)}<C$.\\
Therefore, we have
\begin{equation}\label{eq:a12 gamma2}
\begin{aligned}
a_{12}&=\int_{\partial D_2}\frac{\partial v_1}{\partial
\nu}=\int_{\Gamma_2}\frac{\partial v_1}{\partial
\nu}+\int_{\partial D_2\setminus\Gamma_2}\frac{\partial
v_1}{\partial \nu}=\int_{\Gamma_2}\frac{\partial v_1}{\partial
\nu}+O(1).
\end{aligned}
\end{equation}
 By the harmonicity of $v_1$ on $B(x^0_{12},\delta)\cap
\widetilde\Omega$ and (\ref{eq:v1 gamma3}), we have
\begin{equation}\label{eq:gamma123}
\begin{aligned}
0&=\int_{\Gamma_1}\frac{\partial v_1}{\partial
\nu}+\int_{\Gamma_2}\frac{\partial v_1}{\partial \nu}+
\int_{\Gamma_3}\frac{\partial v_1}{\partial
\nu}=\int_{\Gamma_1}\frac{\partial v_1}{\partial
\nu}+\int_{\Gamma_2}\frac{\partial v_1}{\partial \nu}+O(1).
\end{aligned}
\end{equation}

Meanwhile, by Green's formula and (\ref{eq:v1 gamma3}), we have
\begin{equation} \label{eq:energy B12}
\begin{aligned}
-\int_{B(x^0_{12},\delta)\cap \widetilde\Omega}|\nabla
v_1|^2&=\int_{\Gamma_1}v_1\frac{\partial v_1}{\partial \nu}+
\int_{\Gamma_2}v_1\frac{\partial v_1}{\partial \nu}+\int_{\Gamma_3}v_1\frac{\partial v_1}{\partial \nu}\\
&=\int_{\Gamma_1}\frac{\partial v_1}{\partial \nu}
      +\int_{\Gamma_3}v_1\frac{\partial v_1}{\partial \nu}=\int_{\Gamma_1}\frac{\partial v_1}{\partial \nu}+O(1)
\end{aligned}
\end{equation}
 Therefore, by combining (\ref{eq:a12 gamma2}), (\ref{eq:gamma123}) and  (\ref{eq:energy
 B12}), we have
$$a_{12}=\int_{B(x^0_{12},\delta)\cap \widetilde\Omega}|\nabla v_1|^2+O(1).$$ Similar to the
energy estimates given in Lemma $1.5$, Lemma $1.6$, and Lemma
$1.7$ in \cite{BLY}, we have
\begin{align*}
\frac{1}{C\sqrt {\varepsilon_{12}}}<\int_{B(x^0_{12},\delta)\cap
\widetilde\Omega}|\nabla
v_1|^2<\frac{C}{\sqrt {\varepsilon_{12}}}, ~~~~~~~~~&\text{for}~ n=2\\
\frac{1}{C}|\ln\varepsilon_{12}|<\int_{B(x^0_{12},\delta)\cap
\widetilde\Omega}|\nabla
v_1|^2<C|\ln\varepsilon_{12}|, ~~~~~&\text{for}~ n=3\\
\frac{1}{C}<\int_{B(x^0_{12},\delta)\cap \widetilde\Omega}|\nabla
v_1|^2<C, ~~~~~~~~~~~~~~&\text{for}~ n\geq 4 .
\end{align*}
Therefore,
\begin{align*}
\frac{1}{C\sqrt {\varepsilon_{12}}}<a_{12}<\frac{C}{\sqrt
{\varepsilon_{12}}},~
~~~~~&\text{for}~ n=2,\\
\frac1C|\ln \varepsilon_{12}|<a_{12}<C|\ln\varepsilon_{12}|,~
~~~&\text{for}~ n=3,\\
\frac1C<a_{12}<C,~~~~~~~~~~&\text{for}~ n\geq 4.
\end{align*}$\hfill\square$

Knowing enough properties of the system of linear equations
(\ref{system:C_1 C_2 C_m inclusions}) from Lemma \ref {lm:ab m
inclusions} and Lemma \ref{order aij} , we have
\begin{prop}\label{prop Ci-Cj}
Let $u_\infty\in H^1(\Omega)$ be the weak solution to equation
(\ref{eq:k +infty m inclusions}) and $C_i$ the value of $u_\infty$
on $D_i$, then for any $1\leq i\neq j\leq m$, there exists a
universal constant $C>0$ depending only
 on $n$,  $\kappa_0$,
$r_0$,  $\|\partial \Omega\|_{C^{2,\alpha}}$, $\{\|\partial
D_i\|_{C^{2,\alpha}}\}$, but independent of $\{\varepsilon_{ij}\}$
such that
\begin{equation}
\begin{aligned}
&|C_i-C_j|\leq C\sqrt {\varepsilon_{ij}}\|\varphi\|_{L^\infty(\partial\Omega)}~~~~~~~~~for \hspace{0.2cm} n=2,\\
&|C_i-C_j|\leq C\frac{1}{|\ln{\varepsilon_{ij}}|}\|\varphi\|_{L^\infty(\partial\Omega)}~~~~~~for \hspace{0.2cm} n=3,\\
&|C_i-C_j|\leq
C\|\varphi\|_{L^\infty(\partial\Omega)}~~~~~~~~~~~~~~for
\hspace{0.2cm} n \geq
           4.\\
 \end{aligned}
\end{equation}
\end{prop}
\emph{Proof:}\hspace{.2cm} By Lemma \ref{lm:ab m inclusions}, we
know that the matrix $-A$ satisfies condition $(A1)$ and $(A2)$,
then applying Proposition \ref{prop:system P} on (\ref{system:C_1
C_2 C_m inclusions}), we have, for any $1\leq i\neq j\leq m$,
$$|C_i-C_j|\leq \frac{C}{a_{ij}}\|\varphi\|_{L^\infty(\partial\Omega)}$$
where C is some constant depending on $n$, $\kappa_0$, $r_0$,
$\|\partial \Omega\|_{C^{2,\alpha}}$, $\{\|\partial
D_i\|_{C^{2,\alpha}}\}$, but independent of
$\{\varepsilon_{ij}\}$. \\
By Lemma \ref{order aij}, we immediately finish the
proof.$\hfill\square$

Now we are ready to complete the proof of Theorem \ref{thm:upbdd m
inclusions}.\\
\emph{Proof of Theorem \ref{thm:upbdd m inclusions}:}\hspace{.2cm}
We prove the estimates in dimension 2, the proof for the higher
dimensional cases is similar.  Without loss of generality, we
assume $i=1$, $j=2$ and $\varepsilon_{12}<\delta$. Now we need to
prove the gradient estimates for $u_\infty$ in the narrow region
between $D_1$ and $D_2$. For simplicity, we assume
$\|\varphi\|_{L^{\infty}(\partial\Omega)}=1$.

By the decomposition formula (\ref{rep:u m inclusions}), we have
\begin{equation}
\nabla u_\infty=(C_1-C_2)\nabla
v_1+C_2(\nabla(v_1+v_2))+\sum_{i=3}^m C_i\nabla v_i+\nabla v_0
\nonumber
\end{equation}

By Lemma \ref{lm grad v_0,i m inclusions}, we have
\begin{equation}\label{eq:v0 v1}
\|\nabla v_1 \|_{L^{\infty}(\widetilde{\Omega}\cap
B(x^0_{12},\delta))}<\frac{C}{\varepsilon_{12}},~~~~\|\nabla
v_0\|_{L^{\infty}(\widetilde{\Omega}\cap B(x^0_{12},\delta))}<C
\end{equation}
where $C$ is some universal constant.\\

For $i=3,\ldots, m$, we have, by Lemma \ref{lm grad v_0,i m
inclusions},
\begin{equation}\label{eq:v3}
\|\nabla v_i\|_{L^\infty(\widetilde{\Omega}\cap
B(x^0_{12},\delta))}<C.
\end{equation}

 Since $v_1+v_2=1$ on both $\partial D_1$ and
$\partial D_2$, similar to the proof of Lemma \ref{lm grad v_0,i m
inclusions}, we can show that
\begin{equation}\label{eq:v1+v2}
\|\nabla (v_1+v_2)\|_{L^{\infty}(\widetilde{\Omega}\cap
B(x^0_{12},\delta))}<C.
\end{equation}
By Proposition \ref{prop Ci-Cj}, (\ref{eq:v0 v1}), (\ref{eq:v3})
and (\ref{eq:v1+v2}), we have
\begin{align*}
\|\nabla u_\infty\|_{L^\infty(\widetilde{\Omega}\cap
B(x^0_{12},\delta))}&\leq |C_1-C_2|\|\nabla v_1\|_{
L^\infty(\widetilde{\Omega}\cap
B(x^0_{12},\delta))}+|C_2|\|\nabla(v_1+v_2)\|_{L^\infty(\widetilde{\Omega}\cap
B(x^0_{12},\delta))}\\
&+\sum_{i=3}^m |C_i|\|\nabla
v_i\|_{L^\infty(\widetilde{\Omega}\cap
B(x^0_{12},\delta))}+\|\nabla
v_0\|_{L^\infty(\widetilde{\Omega}\cap B(x^0_{12},\delta))}\\
&\leq C\sqrt{\varepsilon_{12}}\frac{1}{\varepsilon_{12}}+C\\
&\leq\frac{C}{\sqrt{\varepsilon_{12}}}.
\end{align*}

As we mentioned in Remark \ref{rem1}, the strict convexity
assumption of the two inclusions can be weakened. In fact, our
proof of Theorem \ref{thm:upbdd m inclusions} applies, with minor
modification, to more general inclusions as below.

In $\mathbb{R}^n$, $n\geq 2$, for two closely spaced inclusions
$D_i$ and $D_j$ which are not necessarily strictly convex, assume
$\partial D_i\cap B(0,r)$ and $\partial D_j\cap B(0,r)$ can be
represented by the graph of $x_1=f(x')+\frac{\varepsilon_{ij}}{2}$
and $x_1=-g(x')-\frac{\varepsilon_{ij}}{2}$, then $f(0')=g(0')=0$,
$\nabla (g+f)(0')=0$. Assume further that
\begin{equation} \label{expan:g-f}
\lambda_1|x'|^{2l}\leq g(x')+f(x')\leq \lambda_2|x'|^{2l}, ~~~
\forall |x'|\leq r/2,
\end{equation}
where $\lambda_2 \geq \lambda_1>0, l\in \mathbb{Z}^+$.

Under the above assumption, let $u_\infty\in H^1(\Omega)$ be the
solution to equation (\ref{eq:k +infty m inclusions}). Then, for
$\varepsilon_{ij}$ sufficiently small, we have
\begin{equation} \label{bd:R^n}
\begin{split}
&\|\nabla u_\infty\|_{L^{\infty}(\widetilde{\Omega}\cap
B(x^0_{ij},\delta))}\leq C\|\varphi\|_{L^{\infty}(\partial
\Omega)}\varepsilon_{ij}^{-\frac{n-1}{2l}}
~~~~~~~~~~~~~\emph{if}~n-1<2l,\\
&\|\nabla u_\infty\|_{L^{\infty}(\widetilde{\Omega}\cap
B(x^0_{ij},\delta))}\leq C\|\varphi\|_{L^{\infty}(\partial
\Omega)}\frac{1}{\varepsilon_{ij}|\ln\varepsilon_{ij}|}
~~~~~~~~\emph{if}~ n-1=2l,\\
&\|\nabla u_\infty\|_{L^{\infty}(\widetilde{\Omega}\cap
B(x^0_{ij},\delta))}\leq C\|\varphi\|_{L^{\infty}(\partial
\Omega)}\frac{1}{\varepsilon_{ij}} ~~~~~~~~~~~~~~~~~\emph{if}~
n-1>2l.
\end{split}
\end{equation}
where $C$  is a constant depending on $n$, $\lambda_1$,
$\lambda_2$, $r_0$, $\|\partial D_i\|_{C^{2,\alpha}}$ and
$\|\partial D_j\|_{C^{2,\alpha}}$, but independent of
$\varepsilon_{ij}$. \\
For the proof, please refer to the corresponding discussion after
the proof of Theorem 0.1-0.2 in \cite{BLY}.

%%%%%%%%%%%%%%%%%%%%%%%%%%%%%%%%%%%%%%%%%%%%%%%%%%%%%%%%%%%%%%%%%%%%%%%
%
%           insulated conductivity problem
%
%%%%%%%%%%%%%%%%%%%%%%%%%%%%%%%%%%%%%%%%%%%%%%%%%%%%%%%%%%%%%%%%%%%%%%%%%%

\section{The insulated conductivity problem}
In this section, we consider the anisotropic insulated
conductivity problem, which is described by Equation (\ref{eq:k
0}). As we mentioned in the introduction, the gradient can only
blow up when two inclusions are close to each other.  In order to
establish the gradient estimates for this problem, we first
consider the local version of the problem, namely Equation
(\ref{eq:k 0 local harmonic}).

 To make the problem easier, we first consider the equation in a strip. In this case, by
using a ``flipping" technique, we derive the gradient
estimates in the strip. \\
Denote, for any integer $l$
 \begin{equation}
\begin{aligned}
&\mathcal{Q}_l:=\{z\in \mathbb{R}^n\big |
(2l-1)\delta<z_1<(2l+1)\delta, |z'|\leq 1\},\\
&\Gamma_l^+:=\{z\in \mathbb{R}^n\big| z_1= (2l+1)\delta~
\text{and}~ |z'|\leq 1\},\\
&\Gamma_l^-:=\{z\in \mathbb{R}^n\big| z_1= (2l-1)\delta~
\text{and}~ |z'|\leq 1\},
\end{aligned}
\nonumber
\end{equation}
and
$$\mathcal{Q}=\{z\in \mathbb{R}^n \big||z_1|\leq
1~ \text{and}~ |z'|\leq 1\}.$$
 We consider the following equation in $\mathcal{Q}_0$
\begin{equation}\label{eq:w Q 0}
\left\{\begin{aligned}
 &\partial_{z_i}\Big(b^{ij}(z)~\partial_{z_j}{w}\Big)=0~~~~~~~~~
           \text{in}~ \mathcal{Q}_0,\\
 &b^{1j}\partial_{z_j} w=0~~~~~~~~~~~~~~~~~~~~~\text{on}~\Gamma^\pm_0. \\
          \end{aligned}
          \right.
\end{equation}
where $(b^{ij})\in
C^{\alpha}(\overline{\mathcal{Q}}_0)(0<\alpha<1)$ is a symmetric
matrix function in $\mathcal{Q}_0$, and there exist constants
$\Lambda_2\ge\lambda_2>0$ such that, for all $\xi\in\mathbb{R}^n$,
$$\|b^{ij}(z)\|_{C^{\alpha}(\overline{\mathcal{Q}}_0)}\leq \Lambda_2, \hspace{1cm} \lambda_2|\xi|^2\leq b^{ij}(z)\xi_i\xi_j,~~\forall
z\in\mathcal{Q}_0, \xi\in \mathbb{R}^n.$$ Then we have
\begin{lem} \label{lem:w Q 0}
Suppose $w\in H^1(\mathcal{Q}_0)\cap L^\infty(\mathcal{Q}_0)$ is a
weak solution of (\ref{eq:w Q 0}), then there exists a constant
$C>0$ depending only on $n, ~\lambda_2,\Lambda_2 $, but
independent of $\delta$, such that
$$\|\nabla
w\|_{L^{\infty}(\mathcal{Q}_0(\frac{1}{2}))}\leq
C\|w\|_{L^\infty(\mathcal{Q}_0)},$$ where
$\mathcal{Q}_0(\frac{1}{2}):=\{z\in \mathbb{R}^n\big| |z_1|\leq
\delta~ \text{and}~ |z'|\leq \frac{1}{2}\}$.
\end{lem}
\emph{Proof}:\hspace{0.2cm} For any integer $l$,
  We construct a new function $\widetilde w$ by ``flipping" $w$
  evenly in each $\mathcal{Q}_l$.  We define
\begin{equation}
\widetilde w(z)=w((-1)^l(z_1-2l\delta),z'), ~~~~\forall z \in
\mathcal{Q}_{l}.
 \nonumber
\end{equation}
Therefore, we have defined
$\widetilde w$ piecewisely in $\mathcal{Q}$. \\
We define the corresponding elliptic coefficients as follows\\
for $\alpha=2,3,\ldots,n$,
\begin{equation}
\widetilde {b}^{\alpha1}(z)=\widetilde{b}^{1\alpha}(z)=
(-1)^l{b}^{1\alpha}((-1)^l(z_1-2l\delta), z'),~~~~\forall z \in
\mathcal{Q}_l.
 \nonumber
\end{equation}
for all other indices
\begin{equation}
\widetilde {b}^{ij}(z)={b}^{ij}((-1)^l(z_1-2l\delta),
z'),~~~~\forall z \in \mathcal{Q}_l. \nonumber
\end{equation}
Under the above definitions of $\widetilde w$ and $\widetilde
b^{ij}$, we can easily check that, for any integer $l$,
\begin{equation}\label{eq:w Q l}
\left\{\begin{aligned}
 &\partial_{z_i}\Big(\widetilde{b}^{ij}(z)~\partial_{z_j}{\widetilde{w}}\Big)=0~~~~~~~~~~
           \text{in}~ \mathcal{Q}_l,\\
 &\widetilde{b}^{1j}\partial_{z_j} \widetilde{w}=0~~~~~~~~~~~~~~~~~~~~~~\text{on}~\Gamma^\pm_l,
          \end{aligned}
          \right.
          \nonumber
\end{equation}
Then for any test function $\psi\in C^{\infty}_0(\mathcal{Q})$,
we have
\begin{equation}
\begin{aligned}
\int_\mathcal{Q}\widetilde{b}^{ij}(z)~\partial_{z_j}{\widetilde{w}}\partial_{z_i}{\psi}&=
\sum_l\int_{\mathcal{Q}_l}\widetilde{b}^{ij}(z)~\partial_{z_j}{\widetilde{w}}\partial_{z_i}{\psi}\\
&=0 ~~~\big(\textrm{by the definition of weak solution})
\end{aligned}
\nonumber
\end{equation}
Therefore $\widetilde w\in H^1(\mathcal{Q})$ satisfies
\begin{equation} \label{eq:w Q}
\partial_{z_j}\big(\widetilde{b}^{ij}(z)~\partial_{z_i}{\widetilde{w}}\big)=0~~~~~
           \text{in}~ \mathcal{Q}.
           \nonumber
\end{equation}

Following exactly from \cite{LN}, we first introduce a new
equation
\begin{equation} \label{condition:w5}
\partial_{z_i}\big(\widetilde{B}^{ij}(z)~\partial_{z_j}{u}\big)=0~~~~~
           \text{in}~ \mathcal{Q} \nonumber
\end{equation}
where
$$\widetilde B^{ij}(z)=\left\{ \begin{array}{ll}
\lim_{z\in \mathcal{Q}_l, ~z\rightarrow((2l-1)\delta, ~0')}\widetilde {b}^{ij}(z) & z\in\mathcal{Q}_l,l>0;\\
\widetilde {b}^{ij}(0) & z\in\mathcal{Q}_0\\
\lim_{z\in \mathcal{Q}_l, ~z\rightarrow((2l+1)\delta, ~0')}\widetilde b^{ij}(z) & z\in\mathcal{Q}_l,l<0;\\
\end{array} \right.
$$
then we define the norm
$$\|F\|_{Y^{s,p}}=\sup_{0<r<1}r^{1-s}(\int \hspace{-2.5ex}-_{r\mathcal{Q}}{|F|^p})^{\frac{1}{p}}$$
Since $b^{ij}(z)\in C^{\alpha}(\overline{\mathcal{Q}}_0)$ ,
$\widetilde{b}^{ij}(z)$ is piecewise $C^\alpha$ continuous in
$\mathcal{Q}$, then we can immediately check that
$$\|\widetilde {b}^{ij}-\widetilde {B}^{ij}\|_{Y^{1+\alpha,2}}<C$$
where $C$ is some constant only depending on $\Lambda_2$.  Using
Proposition 4.1 in \cite{LN}, we have
$$\|\nabla \widetilde w\|_{L^\infty(\frac{1}{2}\mathcal{Q})}\leq
C\|\widetilde w\|_{L^2(\mathcal{Q})}\leq C\|\widetilde
w\|_{L^\infty(\mathcal{Q})},$$ Then by the definition of
$\widetilde w$, we have
$$\|\nabla
w\|_{L^{\infty}(\mathcal{Q}_0(\frac{1}{2}))}\leq
C\|w\|_{L^\infty(\mathcal{Q}_0)}$$
where $C>0$ depends on $n$,
$\lambda_2$, $\Lambda_2$, but is
independent of $\delta$.  $\hfill\square$ \\

Since $D_1$ and $D_2$ are strictly convex domains, we can write
$\mathcal{O}(r)$, which is defined by (\ref{eq:D1D2}), as follows
 $$\mathcal{O}(r)=\{x\in \mathbb{R}^n\big|
 -g(x')-\varepsilon/2<x_1<f(x')+\varepsilon/2,~ |x'|<r\}$$
 With the side boundary $\Gamma_+$ and $\Gamma_-$ as
 $$\Gamma_+=\{x\in \mathbb{R}^n\big|
 x_1=f(x')+\varepsilon/2, |x'|<r\},~\Gamma_-=\{x\in \mathbb{R}^n\big|
 x_1=-g(x')-\varepsilon/2, |x'|<r\}$$
 where $f(x')$ and $g(x')$ are strictly convex
functions, moreover they satisfy $$f(0')=g(0')=0, ~\nabla
f(0')=\nabla g(0')=0.$$

Under the above notation, we prove Theorem \ref{thm:upbdd k 0}:\\
 \emph{Proof of  Theorem \ref{thm:upbdd
k 0}}:\hspace{.2cm} Fix one point
$(0,x'_0)\in\mathcal{O}(\frac{r}{2})$ and let
$\delta=\sqrt{f(x'_0)+g(x'_0)+\varepsilon}$, since $f(x')$ and
$g(x')$ are strictly convex, then there exists a universal
constant $C$ depending only on $\|\partial D_1\|_{C^{2,\alpha}}$
and $\|\partial D_2\|_{C^{2,\alpha}}$ such that
\begin{equation}\label{eq:x'0}
\frac{1}{C}\sqrt{|x'_0|^2+\varepsilon}<\delta<C\sqrt{|x'_0|^2+\varepsilon}.
\end{equation}

We shift the origin to $(0, x'_0)$ and rescale the coordinates
with $\delta$, then the new coordinates $y=(y_1, y')$ can be
written as follows
\begin{equation}
\left\{\begin{aligned} &y_1=x_1/\delta,\\
 &y'=(x'-x'_0)/\delta.
\end{aligned}
\right.
\end{equation}

Let
$$v(y)=u_0(\delta y_1,x'_0+\delta y'), \hspace{1cm}
\widetilde{a}^{ij}(y)=a^{ij}(\delta y_1,x'_0+\delta y').$$

 Denote
 $$\widetilde{\mathcal{O}}(\widetilde{r}):=\{y\in \mathbb{R}^n\big| -\frac{\varepsilon}{2}-g(x'_0+\delta y')< \delta y_1<
\frac{\varepsilon}{2}+f(x'_0+\delta y'),~|y'|<\widetilde{r}\}$$
With its side boundary
 $$\widetilde{\Gamma}_+:=\{y\in \mathbb{R}^n\big|\delta y_1=\frac{\varepsilon}{2}+f(y'_0+\delta
 y'),
 ~|y'|<\widetilde{r}\}$$
$$\widetilde{\Gamma}_-:=\{y\in \mathbb{R}^n\big|\delta y_1=-\frac{\varepsilon}{2}-g(y'_0+\delta
 y'),
 ~|y'|<\widetilde{r}\}.$$
By (\ref{eq:x'0}), we can find some universal constant
$\widetilde{r}$
 depending only on $\partial D_1$ and $\partial D_2$, such that $\widetilde{\mathcal{O}}(\widetilde{r})$
 is in the image of $\mathcal{O}(r)$ under the above transform.
Thus we have
\begin{equation}
\left\{\begin{aligned} & \partial_{y_i}(\widetilde{a}^{ij}\partial_{y_j}v(y))=0  ~~~~~~\textrm{in}~~~\widetilde{\mathcal{O}}(\widetilde{r}) ,\\
&\widetilde{a}^{ij}\partial_{y_j} v\nu_i=0
~~~~~~~~~~~~~\textrm{on}~~~
\widetilde{\Gamma}_+\cup\widetilde{\Gamma}_-.
\end{aligned}
\right.
\end{equation}
where the coefficients $\widetilde{a}^{ij}$ satisfy, for some
universal constant $C$,
$$\|\widetilde{a}^{ij}\|_{C^\alpha(\widetilde{\mathcal{O}}(\widetilde{r}))}\leq
C\|a^{ij}\|_{C^\alpha(\mathcal{O}(r))}\leq C\Lambda_1,
\hspace{.1cm} \lambda_1|\xi|^2\leq
\widetilde{a}^{ij}(y)\xi_i\xi_j~(\forall
y\in\widetilde{\mathcal{O}}(\widetilde
{r}),~\forall\xi\in\mathbb{R}^n).$$
 Next we construct a map
$\Phi:\widetilde{\mathcal{O}}(\widetilde{r})\longmapsto\mathcal{Q}_0$,
$\Phi(y)=z$ with
\begin{equation}
\left\{\begin{aligned} &z_1=2\delta\frac{\delta y_1+g(x'_0+\delta y')+\varepsilon/2}{f(x'_0+\delta y')+g(x'_0+\delta y')+\varepsilon}-\delta,\\
 &z'= \frac{y'}{\widetilde{r}}.
\end{aligned}
\right.
\end{equation}
It can be verified directly that this map is a diffeomorphism from
$\widetilde{\mathcal{O}}(\widetilde{r})$ to $\mathcal{Q}_0$.

Let
$$w(z)=v(\Phi^{-1}(z))$$ Then from the definition of weak
solution, we know that $w(z)$ satisfies the following equation
\begin{equation}
\left\{\begin{aligned} &\partial_{z_i}\big(b^{ij}(z)\partial_{z_j}w(z)\big)=0  ~~~~~~~\textrm{in}~\mathcal{Q}_0 ,\\
&b^{1j}(z)\partial_{z_j}w(z)=0~~~~~~\textrm{on}~\Gamma_0^+\cup\Gamma_0^-.\\
\end{aligned}
\right.
\end{equation}
where $$\big(b^{ij}(z) \big)=\frac{(\partial_y
z)\big(\widetilde{a}^{ij}(y)\big)(\partial_y
z)^{t}}{|\det(\partial_y z)|}$$

Therefore, we have transferred the original problem into Equation
(\ref{eq:w Q 0}).

 In order to use Lemma \ref{lem:w Q 0}, we have to check that $b^{ij}(z)$ is strictly
elliptic and $\|b^{ij}\|_{C^\alpha(\overline{\mathcal{Q}}_0)}$ is
bounded by some universal constant. First we show that there
exists a universal constant $\lambda_2$ such that
\begin{equation}\label{eq:ellipticity}
\xi^t\big(b^{ij}(z) \big)\xi\ge \lambda_2|\xi|^2
~~~~\forall\xi\in\mathbb{R}^n,~~\forall z\in \mathcal{Q}_0
\end{equation}

Notice that the eigenvalues of $(\partial_y z)$ are
$\frac{1}{\widetilde{r}}$ with multiplicity $n-1$ and
$\partial_{y_1} z_1$. By (\ref{eq:x'0}), we can prove that
\begin{equation}\label{eq:y1z1}
\frac{1}{C}<|\partial_{y_1} z_1|=\partial_{y_1}
z_1=\frac{2\delta^2}{f(x'_0+\delta y')+g(x'_0+\delta
y')+\varepsilon}<C
\end{equation}
where $C$ is some universal constant.

Based on (\ref{eq:y1z1}),we have
$$\xi^t\big(b^{ij}(z) \big)\xi=\xi^t(\partial_y
z)\frac{\big(\widetilde{a}^{ij}(y)\big)}{|\det(\partial_y
z)|}(\partial_y z)^{t}\xi>\lambda_2|\xi|^2, ~~\forall \xi \in
\mathbb{R}^n$$ where $\lambda_2>0$ is some universal constant

  The
boundedness of $\|b^{ij}\|_{C^\alpha(\overline{\mathcal{Q}}_0)}$
can be checked
similarly.\\
 Now applying Lemma \ref{lem:w Q 0}, we have $$\|\nabla
w\|_{L^{\infty}(\mathcal{Q}_0(\frac{1}{2}))}\leq
C\|w\|_{L^\infty(\mathcal{Q}_0)}$$ Tracing back to $u_0$ through
the transforms, we have, for any point
$x\in\mathcal{O}(\frac{r}{2})$,
$$|\nabla u_0(x)|\leq \frac{C\|u_0\|_{L^\infty(\mathcal{O}(r))}}{\delta}\leq \frac{C\|u_0\|_{L^\infty(\mathcal{O}(r))}}{\sqrt{|x'|^2+\varepsilon}}.$$
$\hfill\square$ \\

%%%%%%%%%%%%%%%%%%%%%%%%%%%%%%%%%%%%%%%%%%%%%%%%%%%%%%%%%%%%%%%%%%%%%%%%
%
%      4. Appendix B
%
%%%%%%%%%%%%%%%%%%%%%%%%%%%%%%%%%%%%%%%%%%%%%%%%%%%%%%%%%%%%%%%%%%%%%%%%

\section{Appendix}
{\large\bf Some elementary results for the insulated conductivity
problem}

\vspace{0.3cm} Assume that in $\mathbb{R}^n$, $\Omega$ and
$\omega$ are  bounded open sets with  $C^{2,\alpha}$ boundaries,
$0<\alpha<1$, satisfying, for some $m<\infty$,
$$\overline\omega=\bigcup_{s=1}^{m}\overline\omega_s\subset\Omega,$$
where $\{\omega_s\}$ are
 connected components of $\omega$.
Clearly $\omega_s$ is open for all $1\le s\le m$.
 Given $\varphi\in
C^2(\partial\Omega)$, the conductivity problem we consider is the
following transmission problem with Dirichlet boundary condition:
\begin{equation} \label{diveq:k finite B}
\left\{ \begin{aligned}
           \partial_{x_j}\Big\{\Big[\big(k a_1^{ij}(x)-a_2^{ij}(x)\big)
\chi_\omega+a_2^{ij}(x)\Big]\partial_{x_i}u_k\Big\}&=0~~~in~\Omega, \\
           u_k=\varphi ~~~~~~~~~~~~~~~~~~~&~~~~~~~~on~\partial\Omega,
          \end{aligned}
\right.
\end{equation}
where $0<k<1$, and $\chi_\omega$ is the characteristic function of
$\omega$.

The $n\times n$ matrixes
$A_1(x):=\big(a_1^{ij}(x)\big)~\text{in}~\omega,~
A_2(x):=\big(a_2^{ij}(x)\big)~\text{in}~\Omega\backslash\overline{\omega}$
are symmetric and $\exists$ a constant $\Lambda\geq\lambda>0$ such
that
$$~\lambda|\xi|^2\leq a_1^{ij}(x)\xi_i\xi_j\leq \Lambda|\xi|^2~(\forall x\in\omega),
~~~~\lambda|\xi|^2\leq a_2^{ij}(x)\xi_i\xi_j\leq
\Lambda|\xi|^2~(\forall x\in\Omega\backslash\omega)$$ for all
$\xi\in\mathbb{R}^n$ and $a_1^{ij}(x)\in C^2(\overline \omega),
~a_2^{ij}(x) \in C^2(\overline \Omega\backslash\omega)$.

\vspace{0.2cm} Equation (\ref{diveq:k finite B}) can be rewritten
in the following form to emphasize the transmission condition on
$\partial\omega$:
\begin{equation} \label{eq:k finite B}
\left\{ \begin{aligned}
           &\partial_{x_j}\Big(a_1^{ij}(x)~\partial_{x_i}{u_k}\Big)=0~~~~~~~~~~~~~~~~~~~~in~\omega,\\
           &\partial_{x_j}\Big(a_2^{ij}(x)~\partial_{x_i}{u_k}\Big)=0~~~~~~~~~~~~~~~~~~~~in~\Omega
           \backslash\overline\omega,\\
           &u_k|_{+}=u_k|_{-},~~~~~~~~~~~~~~~~~~~~~~~~~~~~~~~on~\partial\omega, \\
           &a_2^{ij}(x)\partial_{x_i}{u_k}\nu_j\big|_{+}=ka_1^{ij}(x)\partial_{x_i}{u_k}\nu_j\big|_{-}~~~on~\partial\omega, \\
           &u_k=\varphi~~~~~~~~~~~~~~~~~~~~~~~~~~~~~~~~~~~~~~~on~\partial\Omega.
          \end{aligned}
\right.
%\nonumber
\end{equation}

It is well known that  equation (\ref{diveq:k finite B}) has a
unique  solution $u_k$ in $H^1(\Omega)$, and the solution $u_k$ is
in  $C^1(\overline{\Omega\backslash\omega})\cap
C^1(\overline\omega)$ and satisfies equation (\ref{eq:k finite
B}). On the other hand, if
 $u_k\in C^1(\overline{\Omega\backslash\omega})\cap
C^1(\overline\omega)$ is a solution of equation (\ref{eq:k finite
B}), then $u_k\in H^1(\Omega)$  satisfies equation (\ref{diveq:k
finite B}).

 For  $k\in (0, 1)$, consider the energy functional
\begin{equation} \label{vm:k finite B}
\begin{aligned}
I_{k}[v]:&
=\frac{k}{2}\int_{\omega}a_1^{ij}(x)\partial_{x_i}{v}\partial_{x_j}{v}
+
\frac{1}{2}\int_{\Omega\backslash\overline{\omega}}a_2^{ij}(x)\partial_{x_i}{v}\partial_{x_j}{v},
\end{aligned}
%\nonumber
\end{equation}
defined on
$$H^1_{\varphi}(\Omega)
:=\{v\in H^1(\Omega)| ~v=\varphi ~~on~ \partial\Omega\}.$$

It is well known that for
  $k\in (0, 1)$, the solution $u_k$ of (\ref{diveq:k finite B})
is the  minimizer of the minimization problem:
$$I_{k}[u_k]=\min_{v\in H^1_{\varphi}(\Omega)}I_{k}[v].$$

For $k=0$, the insulated conducting problem is:
\begin{equation} \label{eq:k 0 B}
\left\{ \begin{aligned}
           &\partial_{x_j}\Big(a_2^{ij}(x)~\partial_{x_i}{u_0}\Big)=0~~~~~~~~~~~~~~~~~~~~in~\Omega
           \backslash\overline\omega,\\
           &a_2^{ij}(x)\partial_{x_i}{u_0}\nu_j\big|_{+}=0~~~~~~~~~~~~~~~~~~~~~~~on~\partial\omega, \\
           &u_0=\varphi~~~~~~~~~~~~~~~~~~~~~~~~~~~~~~~~~~~~~~~on~\partial\Omega,\\
           &\partial_{x_j}\Big(a_1^{ij}(x)~\partial_{x_i}{u_0}\Big)=0~~~~~~~~~~~~~~~~~~~~in~\omega,\\
           &u_0|_{+}=u_0|_{-},~~~~~~~~~~~~~~~~~~~~~~~~~~~~~~~~~on~\partial\omega. \\
          \end{aligned}
\right.
%\nonumber
\end{equation}

Equation  (\ref{eq:k 0 B}) has a unique solution $u_0\in
H^1(\Omega)$, which can be solved in $\Omega\setminus
\overline\omega$ by the first three lines in   (\ref{eq:k 0 B}),
and then, with $u_0|_{\partial \omega}$, be solved in $\omega$
using the fourth line in
  (\ref{eq:k 0 B}).
It is well known that $u_0\in C^1(\overline \Omega\setminus
\omega)\cap C^1(\overline \omega)$.

Define the energy functional
\begin{equation} \label{vm:k 0 B}
I_{0}[v]:=
\frac{1}{2}\int_{\Omega\backslash\overline{\omega}}a_2^{ij}(x)\partial_{x_i}{v}\partial_{x_j}{v},
%\nonumber
\end{equation}
where $v$ belongs to the set
$$
\mathcal{A}_0:=\big\{v\in H^1(\Omega\setminus \overline \omega)
\big| \ v=\varphi \ \mbox{on}\ \partial \Omega\}.
$$

It is well known that there is a unique $v_0\in \mathcal{A}_0$
which is the minimizer to the minimization problem:
$$I_{0}[v_0]=\min_{v\in \mathcal{A}_0}I_{0}[v].$$
Moreover, $v_0=u_0$ a.e. in $\Omega\setminus \overline \omega$,
 where
$u_0$ is the solution of
  (\ref{eq:k 0 B}).

Now, we give the relationship between $u_k$ and $u_0$.
\begin{thm} \label{thm:weak limit B}
For $0<k<1$, let $u_k$ and $u_0$ in $H^1(\Omega)$ be the solutions
of equations (\ref{eq:k finite B}) and (\ref{eq:k 0 B}),
respectively. Then
\begin{equation}
u_k\rightharpoonup u_0 ~~\text{in}~H^1(\Omega), ~~~\text{as}~
k\rightarrow 0, \label{1}
\end{equation}
 and, consequently,
\begin{equation}
\lim_{k\rightarrow 0}I_{k}[u_k] =I_{0}[u_0]. \label{2}
\end{equation}
\end{thm}

\emph{Proof}:\hspace{0.2cm} We will first show that
\begin{equation}
\sup_{ 0<k<1} \|\nabla u_k\|_{ L^2(\Omega)}<\infty. \label{bound}
\end{equation}

Since $u_k$ is the minimizer of $I_{k}$ in $H^1_{\varphi}(\Omega)$
and $v_0:=u_0|_{\Omega\setminus\overline\omega}$
 is the minimizer of $I_{0}$ in $\mathcal{A}_0$, we have
\begin{equation}
\begin{aligned}
\frac{\lambda k}{2}\|\nabla u_k\|_{L^2(\omega)}+I_{0}[v_0]
&\leq\frac{k}{2}\int_{\omega}a_1^{ij}(x)\partial_{x_i}{u_k}
\partial_{x_j}{u_k}+I_{0}[v_0]\\
&\leq
\frac{k}{2}\int_{\omega}a_1^{ij}(x)\partial_{x_i}{u_k}\partial_{x_j}{u_k}
+I_{0}[u_k|_{\Omega\setminus\overline\omega}]
=I_{k}[u_k]\\
&\leq
I_{k}[u_0]=\frac{k}{2}\int_{\omega}a_1^{ij}(x)\partial_{x_i}{u_0}
\partial_{x_j}{u_0}+I_{0}[v_0],\\
&\leq \frac{\Lambda k}{2}\|\nabla u_0\|_{L^2(\omega)}+I_{0}[v_0] .
\end{aligned}
\nonumber
\end{equation}
Thus
$$\sup_{ 0<k<1}
\|\nabla u_k\|_{L^2(\omega)} <\infty.
$$

On the other hand,
$$\frac{\lambda}{2}\|\nabla
u_k\|_{L^2(\Omega\backslash\overline\omega)}\leq I_{k}[u_k]\leq
I_{k}[u_0]\leq \frac{\Lambda}{2}\|\nabla u_0\|_{L^2(\Omega)}.$$

Estimate (\ref{bound}) follows from the above.

Since $u_k=\varphi$ on $\partial\Omega$, we derive from
(\ref{bound}) that $\sup_{0<k<1}\|u_k\|_{H^1(\Omega)}<\infty$. Let
$u_k\rightharpoonup u^*_0
\hspace{0.2cm}in\hspace{0.2cm}H^1_{\varphi}(\Omega)$ along a
subsequence of  $k\rightarrow 0$
(still denoted as $k\to 0$).\vspace{0.2cm}\\
We will show that
 $u^*_0$ is  a solution of equation (\ref{eq:k 0
B}). Therefore, $u^*_0= u_0$.

\vspace{0.1cm}We  only need to establish the following three
properties:
\begin{align}
\label{condition:u a2}\partial_{x_j}\Big(a_2^{ij}(x)~\partial_{x_i}{u^*_0}\Big)&=0~~~~~~~~in~\Omega\backslash\overline\omega,\\
\label{condition:u a1}\partial_{x_j}\Big(a_1^{ij}(x)~\partial_{x_i}{u^*_0}\Big)&=0~~~~~~~~in~\omega,\\
\label{condition:0} u^*_0\in  C^1(\Omega\setminus \omega), \qquad
a_2^{ij}(x)\partial_{x_i}{u^*_0}\nu_j\big|_{+}&=0
~~~~~~~~on~\partial\omega.
\end{align}
(i) For  $k\in (0, 1)$, we see from equation (\ref{diveq:k finite
B}) that
$$
\partial_{x_j}\Big(a_2^{ij}(x)~\partial_{x_i}{u_k}\Big)=0,\qquad
\mbox{in}\ \Omega\setminus \overline \omega,
$$
$$
\partial_{x_j}\Big(a_1^{ij}(x)~\partial_{x_i}{u_k}\Big)=0,\qquad
\mbox{in}\ \omega.
$$
Since $u_k$ converges to $u^*_0$ weakly in $H^1(\Omega)$,
(\ref{condition:u a2}) and (\ref{condition:u a1}) follow from the
above.

(ii) For any $w\in {\cal A}_0$, we extend it to $\tilde w\in
H^1_\varphi(\Omega)$ (i.e. $\tilde w=w$ in $\Omega\setminus
\overline \omega$). By the minimality of $u_k$,
$$
I_k(u_k)\le I_k(\tilde w).
$$
Sending $k$ to $0$ leads to
$$
I_0(u^*_0|_{ \Omega\setminus \omega})\le I_0(w).
$$
Thus $u^*_0=u_0$ a.e. in $\Omega\setminus \omega$.
(\ref{condition:0}) follows.

We have proved (\ref{1}). Theorem \ref{thm:weak limit B} is
established.


\begin{thebibliography}{99}

\bibitem{AKLLL} H. Ammari, H. Kang, H. Lee,
 J. Lee and M. Lim, {\it Optimal estimates for the Electrical Field in Two
 Dimensions}, J. Math. Pures Appl. 88 (2007), 307-324.

\bibitem{AKLLZ}  H. Ammari, H. Kang,  H. Lee,
 M. Lim and H. Zribi,
{\it Decomposition Theorems and Fine Estimates for
Electrical Fields in the Presence of Closely Located},
preprint.

\bibitem{AKL} H. Ammari, H. Kang and M. Lim, {\it Gradient estimates for solutions to the conductivity problem},
Math. Ann. 332 (2005), 277-286.

\bibitem{BASL}  I. Babuska, B. Anderson, P.J. Smith,
and K. Levin, {\it Damage analysis of fiber composites. I.
Statistical analysis on fiber scale}, Comput. Methods Appl. Mech.
Engrg. 172 (1999),  27-77.

\bibitem{BLY} E.S. Bao, Y.Y. Li and B. Yin, {\it Gradient estimates for the
perfect conductivity problem}, Arch. Rational Mech. and Anal. 193
(2009), 195-226.

\bibitem{BV} E. Bonnetier and M. Vogelius, {\it An elliptic regularity result for a
composite medium with ``touching" fibers of circular
cross-section}, SIAM J. Math. Anal. 31 (2000),  651-677.

\bibitem{BC} B. Budiansky and G.F. Carrier, {\it High shear stresses in stiff fiber composites}, J. App.
Mech. 51 (1984), 733-735.

\bibitem{F} G.B. Folland, {\it Introduction to Partial Differential Equations},
Princeton University Press, Princeton, NJ, 1976.

\bibitem{GT} D.Gilbarg and N. Trudinger, {\it Elliptic Partial
Differential Equations of Second Order}, Springer, 1998.

\bibitem{GV} R. Horn and C. Johnson, {\it Matrix Analysis}, Cambridge University Press,
1985.

\bibitem{K} J.B. Keller, {\it Stresses in narrow regions}, Trans. ASME J. Appl.
Mech. 60 (1993),  1054-1056.

\bibitem{KSW} D. Kapanadze and B.W. Schulze, {\it Boundary-contact problems for domains with edge singularities},
  J. Differential Equations  234  (2007),  26-53.

\bibitem{LN}  Y.Y. Li and L. Nirenberg, {\it Estimates for elliptic
system from composite material}, Comm. Pure Appl. Math. 56 (2003),
892-925.

\bibitem{LV}  Y.Y. Li and M. Vogelius, {\it Gradient estimates for
solution to divergence form elliptic equation with discontinuous
coefficients}, Arch. Rational Mech. Anal. 153 (2000), 91-151.

\bibitem{LY} M. Lim and K. Yun,  {\it Blow-up of electric fields
between closely spaced spherical perfect conductors},
Communications in PDE, to appear.

\bibitem{MOV} J. Mateu, J. Orobitg
 and J. Verdera, {\it Extra cancellation of even Calderon-Zygmund operators
and quasiconformal mappings},  J. Math. Pures Appl.  91  (2009),
 402-431.

\bibitem{V}  V.V. Voyevodin, {\it Linear Algebra}, Mir Publishers,
Moscow,  1983.

\bibitem{M} X. Markenscoff, {\it Stress amplification in
vanishingly small geometries}, Computational Mechanics 19 (1996),
77-83.

\bibitem{Yu}  K. Yun, {\it Estimates for electric fields blown up
between closely adjacent conductors with arbitrary shape}, SIAM J.
on Applied Math. 67 (2007), 714-730.

\bibitem{Yu2} K. Yun, {\it Optimal bound on high stresses occurring between
stiff fibers with arbitrary shaped cross-sections}, J. Math. Anal.
Appl. 350 (2009),  306--312.
\end{thebibliography}
\end{document}